\documentclass[12pt]{article}
\usepackage{amssymb,latexsym, amsmath, amscd, array, graphicx}
\makeatletter
\date{}

\newtheorem{theorem}{Theorem}[section]

\newtheorem{proposition}[theorem]{Proposition}
\newtheorem{definition}[theorem]{Definition}

\newtheorem{corollary}[theorem]{Corollary}

\newtheorem{problem}[theorem]{Problem}
\newtheorem{conjecture}[theorem]{Conjecture}
\newcommand{\edim}{{\rm e}$-${\dim}}

\newcommand{\z}{{\Bbb Z}}

\newcommand{\q}{{\Bbb Q}}
\newcommand{\re}{{\Bbb R}}
\newcommand{\N}{{\Bbb N}}

\newcommand{\st}{{\rm st}}

\newcommand{\invlim}{{\rm invlim}}
\newcommand{\no}{$\rm  {N\ddot{o}beling}$ }

\newcommand{\lo}{\longrightarrow}

\newcommand{\black}{{\blacksquare}}
\newcommand{\tor}{{\rm Tor}}
\newcommand{\diam}{{\rm diam}}

\begin{document}

\title{On compacta not admitting a stable intersection in   $\re^n$}

\author{  Michael Levin}

\maketitle
\begin{abstract} Compacta $X$ and $Y$ are said to admit a stable  intersection in ${\re^n}$ if
 there are  maps $f : X \lo \re^n$ and $g : Y \lo \re^n$ such that 
 for every sufficiently close continuous  approximations 
 $f' : X \lo \re^n$ and $g' : Y \lo \re^n$ of $f$ and $g$ 
 we have  $f'(X)\cap g'(Y)\neq\emptyset$. The well-known conjecture
 asserting  that $X$ and $Y$ do not admit a stable intersection
  in $\re^n$ if and only if 
 $\dim X \times Y \leq n-1$ was confirmed  in many cases. In this paper
 we prove this conjecture in  all the remaining cases except the case $\dim X =\dim Y =3$, $\dim X \times Y=4$ 
 and $n=5$ which still remains open.
\bigskip
\\
{\bf Math. Subj. Class.:}  55M10 (54F45 55N45)
\end{abstract}
\begin{section}{Introduction}
All the spaces are assumed to be separable  metrizable. A map means
a continuous map and a compactum means a compact metric space. 
Compacta $X$ and $Y$ are said to admit a stable  intersection in ${\re^n}$ if
 there are  maps $f : X \lo \re^n$ and $g : Y \lo \re^n$ such that 
 for every sufficiently close continuous  approximations  
  $f' : X \lo \re^n$ and $g' : Y \lo \re^n$ of $f$ and $g$
 we have  $f'(X)\cap g'(Y)\neq\emptyset$.
 The modern 
 research on this subject was initiated
 by a work of D. McCullough
and L. Rubin \cite{m-rubin} who improved the classical 
\no-Pontrjagin theorem by showing that every map
of an $n$-dimensional Boltyanskii compactum $X$ to $\re^{2n}$ 
can be arbitrarily closely approximated
by an embedding. Recall that a compactum $X$ is called a Boltyanskii
compactum if  $\dim X^2 <2\dim X$. One can easily observe 
that any map from a compactum $X$ to $\re^{2n}$ can be approximated 
by embeddings if and only if $X$ does not admit a stable
intersection with itself  in $\re^{2n}$.
This motivated the following  well-known conjecture.
\begin{conjecture}
\label{conjecture}
Compacta   $X$ and $Y$ do not admit a stable  intersection
 in $\re^n$ if and only
if $\dim X \times Y \leq n-1$.
\end{conjecture}

An extensive work on this conjecture  culminated
in the following results.

\begin{theorem}{\rm ({\bf Dranishnikov-West} \cite{dranish-west})}
\label{stable}
Let $X$ and $Y$ compacta such that $\dim X \times Y \geq n$. Then
$X$ and $Y$ admit a stable  intersection in $\re^n$.

\end{theorem}
\begin{theorem}
{\rm ({\bf Dranishnikov-Repovs-Schepin \cite{d-r-s}, Torunczyk-Spiez \cite{t-s}})}
\label{complementary}
Let $X$ and $Y$ be compacta such that $2\dim X + \dim Y\leq 2n-2$ and
$\dim X \times Y \leq n-1$.
Then $X$ and $Y$ do not admit a stable intersection in $\re^n$.
\end{theorem}

\begin{theorem}
\label{co-dim-3}
{\rm ({\bf Dranishnikov }\cite{dranishnikov-intersection})}
Let $X$ and $Y$ be compacta such that $\dim X \leq n-3$, $\dim Y\leq n-3$
and $\dim X \times Y \leq n-1$. Then $X$ and $Y$ do not admit 
a stable intersection in $\re^n$.
\end{theorem}
The goal of this paper is to settle all the remaining open cases
of Conjecture \ref{conjecture}  except only one which still remains open.

\begin{theorem}
\label{main}
Conjecture \ref{conjecture} holds in all the cases 
except the following  one
which still remains open: $\dim X =\dim Y =3$, $\dim X \times Y=4$
and $n=5$.
\end{theorem}

The paper is organized as follows: basics 
of Cohomological Dimension are presented in 
Section \ref{cohomological-dimension}, 
 maps to CW-complex are discussed 
in Section \ref{cw-complexes} and applied  
to obtain
a factorization theorem,
Theorem \ref{main} is proved in Section \ref{main-result}
and a few remarks related to the paper's results are given
in the last section.

\end{section}
\begin{section}{Cohomological Dimension}
\label{cohomological-dimension}
Let us  review basic facts  of cohomological dimension.
By   cohomology we always mean the Cech
cohomology. Let $G$ be an abelian group. The  cohomological dimension
 $\dim_GX$
of a space $X$ with respect to the coefficient group $G$ does not exceed $n$, $\dim _G X \leq n$ if $H^{n+1}(X,A;G)=0$ for every closed $A
\subset X$. We note that this condition implies that
$H^{n+k}(X,A;G)=0$ for all $k\ge 1$ \cite{Ku},\cite{DrArxiv}.
Thus, $\dim _G X =$ the smallest integer $n\geq 0$ satisfying
$\dim _G X \leq n$ and $\dim _G X = \infty $ if such an integer
does not exist. Clearly, $\dim_G X \leq \dim_{\z}X\le\dim X$.

 \begin{theorem}
 \label{alexandrov}
 {\rm ({\bf Alexandroff})}
$\dim X=\dim_\z X$ if $X$ is a finite dimensional space.
\end{theorem}

Let $\mathcal P$ denote the set of all primes. The {\em  Bockstein basis} is the collection of groups
$\sigma=\{ \q, \z_p , \z_{p^\infty}, \z_{(p)} \mid p\in\mathcal P
\}$ where $\z_p =\z/p\z$ is the $p$-cyclic group,
$\z_{p^\infty}={\rm dirlim} \z_{p^k}$  is the $p$-adic circle, and
$\z_{(p)}=\{ m/n \mid n$ is not divisible by $p \}\subset\q$  is the $p$-localization of integers.
\\\\
 The Bockstein basis   of an abelian group $G$ is the collection
$\sigma(G) \subset \sigma$ determined by the rule:

$\z_{(p)} \in \sigma(G)$ if $G/\tor G$ is not divisible by $p$;

$\z_p \in \sigma(G)$ if $p$-$\tor G$ is not divisible by $p$;

$\z_{p^\infty} \in \sigma(G)$ if $p$-$\tor G\neq 0$ is  divisible by $p$;

$\q \in \sigma(G)$ if $G/\tor G\neq 0$ is  divisible by all $p$.
\\
\\
Thus  $\sigma(\z)=\{\z_{(p)}\mid p\in\mathcal P\}$.

\begin{theorem}
\label{bockstein-theorem}
{\rm ({\bf  Bockstein Theorem})} For a compactum $X$,
$$\dim_G X =\sup \{ \dim_H X : H \in \sigma(G)\}.$$
\end{theorem}
${}$\\
Suggested by the Bockstein inequalities we say that
a function $D:\sigma\lo \N \cup \{0, \infty\}$ is   a {\em
$p$-regular} if
$
D(\z_{(p)})=D(\z_p)=D(\z_{p^{\infty}})=D(\q)$ and it  is   {\em
$p$-singular} if
$ D(\z_{(p)})=\max\{ D(\q),D(\z_{p^{\infty}})+1\}.$
A $p$-singular function $D$ is called $p^+$-singular if
$D(\z_{p^{\infty}})=D(\z_p)$ and it is
called {\em $p^-$-singular} if
$D(\z_{p^{\infty}})=D(\z_p)-1$. A function
$D:\sigma\lo \N \cup \{0, \infty\}$ is called a {\em dimension
type} if for every prime $p$ it is either $p$-regular or
$p^{\pm}$-singular. 
Thus, the values of $D(F)$ for the Bockstein fields $F\in\{\z_p,\q\}$ 
together with $p$-singularity types of $D$
determine the value $D(G)$ for all groups in $\sigma$.
For a dimension type $D$ denote $\dim D =\sup \{ D(G) : G \in \sigma \}$.
Note that
  $D(G)\geq 1 $
 for every $G \in \sigma$ if $\dim D >0$. 

\begin{theorem}
\label{bockstein-inequalities}
{\rm ({\bf Bockstein Inequalities} \cite{Ku}, \cite{DrArxiv})}
For every space $X$ the function $d_X:\sigma\lo \N \cup \{0,
\infty\}$ defined as $d_X(G)=\dim_GX$ is a dimension type. 
\end{theorem}
The function
$d_X$ is called the {\em dimension type of $X$}.

\begin{theorem}
\label{dranishnikov-realization-theorem}
{\rm ({\bf Dranishnikov  Realization Theorem} \cite{DrUsp},\cite{DrPacific})}
For every   dimension type $D$ there is a compactum $X$
with $d_X=D$ and $\dim X=\dim D$.
\end{theorem}

\begin{theorem}
{\rm({\bf Olszewski Completion Theorem}\cite{Ol})}
\label{completion}
For every space $X$ there is a complete space $X'$ such that
$X \subset X'$ and $d_X =d_{X'}$.
\end{theorem}

Let $D$ be a dimension type. We will use the abbreviations
 $D(0)=D(\q)$, $D(p)=D(\z_p)$. Additionally, if $D(p)=n\in\N$ we will write $D(p)=n^+$ if $D$ is $p^+$-regular
and $D(p)=n^-$ if it is $p^-$-regular. For a $p$-regular $D$  we leave it without decoration: $D(p)=n$.
Thus, any sequence of decorated numbers $D(p)\in\N$, where $p\in\mathcal P\cup\{0\}$  define a unique dimension type.
There is a natural order on decorated numbers $$\dots<n^-<n<n^+<(n+1)^-< \dots\ .$$
Note that the inequality of dimension types $D\le D'$ as functions on $\sigma$ is equivalent to the family of inequalities
$D(p)\le D'(p)$ for the above order for all $p\in\mathcal P\cup\{0\}$.
Also note that $0$  has no decoration,  $1$ does not have  the '-' decoration
and 
 $D(0)=D(\q)$ has no decoration. 

Let $\epsilon$ be a decoration. We  define 
the reversed decoration $-\epsilon$ and the commutative product  
of decorations $\otimes$ as follows:

$$-(-)=+, -(+)=-,  -({\rm no \ decoration})={\rm no \ decoration}$$
$$\epsilon\otimes ({\rm no \ decoration})=\epsilon,\ \  \ \ \epsilon\otimes \epsilon=\epsilon, \ \ \  \text{and}\ \ \
+\otimes -=-.$$

For dimension types $D_1$ and $D_2$  we define
  the dimension types $D_1 \boxplus D_2$  and $D_1\oplus D_2$
 as follows:
 if $D_1(p)=n^{\epsilon_1}$ and $D_2(p)=m^{\epsilon_2}$
 where $\epsilon_i$ is a decoration then 
 $$(D_1\boxplus D_2)(p)=(n+m)^{\epsilon_1\otimes\epsilon_2}$$
$$(D_1\oplus D_2)(p)=(n+m)^{-((-\epsilon_1)\otimes(-\epsilon_2))}.$$

For an integer  $n\geq 0$ we denote by $n$ the dimension type which sends every
 $G\in \sigma$ to $n$  and
  for a dimension type $D$ we denote  by $D +n$  the dimension type
  which is  the ordinary sum
 of $D$ and $n$ as functions. Note that $D +n$ 
 preserves the decorations of $D$. Also note that $d_{\re^n}=n$.
 
 The operations  $D_1 \boxplus D_2$ and $D_1 \oplus D_2$
  are   motivated by
the following properties.
\begin{theorem}
\label{bockstein-product-theorem}
{\rm(\bf {Bockstein Product Theorem} \cite{DrArxiv},\cite{Sch}, \cite{DyA})}
For  any two compacta  $X$ and $Y$ 

$$d_{X \times Y} =d_X \boxplus d_Y.$$
\end{theorem}

 \begin{theorem}
 \label {union}
 {\rm ({\bf Dydak Union Theorem}\cite{Dy},\cite{first-exotic})}
 Let $X$ be a compactum and $D_1$ and $D_2$  dimension types and let 
 $X=A \cup B$ be  a decomposition with $d_A \leq D_1$ and
 $d_B \leq D_2$.  Then $d_X\leq D_1 \oplus D_2 +1$.
\end{theorem}

\begin{theorem}
\label{decomposition}
{\rm ({\bf Dranishnikov  Decomposition  Theorem} \cite{DrPacific}, \cite{first-exotic})}
  Let $X$ be a finite dimensional compactum 
  and $D_1$ and $D_2$ dimension types such that   $d_X \leq D_1 \oplus D_2 +1$.
 Then there is a decomposition  $X=A \cup B$ such that
 $d_A \leq D_1$ and   $d_B \leq D_2$.
\end{theorem} 
 For a dimension type $D$  and $n \geq \dim D$ we define
  the dimension type $n+1 \ominus D$ by: $$(n+1\ominus D)(p)=(n+1-m)^{-\epsilon}
  \ \ {\rm if}\ \ D(p)=m^\epsilon.$$
  	Note that that $n+1\ominus D$ is indeed  a dimension type and if $\dim D >0$
	then 
	$\dim (n+1\ominus D) \leq n $ and $n+1-(n+1\ominus D)=D$.
  One can also easily verify 
   the  following properties:
  \begin{proposition}
  \label{type-minus}
  Let $D$ be a dimension type and $n=\dim D$. Then
  
  (i) $D \boxplus (n+1 \ominus D)  \leq n+1$. Moreover, for 
  a dimension type $D'$ the condition  $\dim (D \boxplus D') \leq n+1$
  imply that $D'\leq n+1\ominus D$;
  
  (ii) $n+1 \leq  D\oplus (n+1\ominus D) $. Moreover,  for  dimension
  types $D'$ and $D''$ the conditions  $D'\leq D$, $D'' \leq n+1\ominus D$
  and $ n+1 \leq D'\oplus D''$ imply  that $D'=D$ and $D''=n+1\ominus D$.
  \end{proposition}
  In particular (ii) of Proposition \ref{type-minus}
    implies Theorem 
    \ref{dranishnikov-realization-theorem}.
     Indeed, let $D$ be
   a dimension    type and  $n=\dim D$. By \ref{decomposition} consider
   a decomposition  $A \cup B=\re^{n+2}$  with
   $d_A \leq D$ derived 
    from the inequality $n+2 \leq D\oplus (n+1\ominus D) +1 $.
  Then,
   by (ii) of \ref{type-minus} and Theorems  \ref{completion} and
   \ref{union}, one can easily conclude that $\re^{n+2}$ contains
  a  compact subset  of dimension type $D$. 
	
\end{section}

\begin{section}{Maps to CW-complexes}
\label{cw-complexes}
The goal of this section is to prove
\begin{theorem} 
\label{factorization} Any map from a finite dimensional compactum $X$ 
to   a finite CW-complex $L$ with $\dim L \leq 3$
 can be arbitrarily closely approximated by
a map  which factors
through a compactum $Z$ with $d_{Z} \leq d_X$ and $\dim Z \leq 3$.
\end{theorem}
 For proving this theorem 
we need a few facts from  Extension Theory.

Cohomological Dimension is characterized by the following basic property:
$\dim_G X \leq n$ if and only  for every closed
$A \subset X$ and a map $f : A \lo K(G,n)$,
$f$ continuously extends over $X$ where $K(G,n)$ is the Eilenberg-MacLane complex
of type $(G,n)$
(we assume that $K(G,0)=G$ with discrete topology and $K(G, \infty)$ is a singleton).
This extension characterization of Cohomological
Dimension gives a rise to  Extension  Theory (more general than
Cohomological Dimension Theory)
 and  the notion of Extension Dimension.
 The {\em extension dimension} of a space $X$ is said
to be dominated by a CW-complex $K$, written $\edim X \leq K$, if
every map $f : A \lo K$ from a closed subset $A$ of $X$
continuously extends over $X$. Thus $\dim_G X \leq n$ is equivalent
to $\edim X \leq K(G,n)$ and $\dim X \leq n$ is equivalent to 
$\edim X \leq S^n$.

The following theorem shows a  close connection between extension
and cohomological dimensions.

\begin{theorem}
\label{dranishnikov-extension-theorem}
{\rm (\bf {Dranishnikov Extension Theorem} \cite{dranishnikov-sbornik},\cite{Dy})}
Let $K$ be a  CW-complex and
  $X$ a  metric space. Denote by $H_*(K)$ the reduced integral homology
  of $K$. Then

(i)
$\dim_{H_n(K)} X \leq n$ for every $n\geq 0$
 if
$\edim X \leq K$;

(ii) $\edim X \leq K$ if $K$ is simply connected, $X$ is finite dimensional
and
$ \dim_{H_n(K)} X \leq n$ for every $n\geq 0$.
\end{theorem}

Let $G$ be an abelian group.  
We always assume that 
a Moore space $M(G,n)$ of  type $(G,n)$ is an $(n-1)$-connected CW-complex.
Theorem \ref{dranishnikov-extension-theorem} implies
that for a finite dimensional compactum $X$ and $n>1$,  $\dim_G X \leq n$
if and only if $\edim X \leq M(G,n)$. We will refer to this property
of $M(G,n)$ as being a classifying space for finite dimensional
compacta of  $\dim_G \leq n$.  We are going to extend this property to 
Moore spaces   $M(G,1)$ for the groups $G$ in the Bockstein 
basis $\sigma$.

We will consider the following standard models of $M(G,1)$, $G \in \sigma$:
\begin{itemize}
\item
$M(\q, 1)$=the infinite telescope  of  a sequence of maps  from
$S^1 \lo S^1$ of all possible non-zero degrees; 
\item
$M(\z_{(p)},1)$= the infinite telescope of  a sequence of maps  
$S^1\lo S^1$ of all possible non-zero degrees not divisible by $p$;
\item
$M(\z_{p^\infty}, 1)$= the infinite telescope of 
an $p$-fold covering  map  $S^1\lo S^1$  with a disk attached
to the first circle of the telescope by the identity map
of the disk boundary;
\item
$M(\z_p,1)$= a disk attached to $S^1$ by an $p$-fold covering map
from the disk boundary to $S^1$.
\end{itemize}
Note that  $M(\q,1)=K(\q,1)$ and $M(\z_{(p)}, 1)=K(\z_{(p)},1)$
and hence $M(\q,1)$ and $M(\z_{(p)}, 1)$ are classifying spaces for
compacta of $\dim_\q \leq 1$ and 
$\dim_{\z_{(p)}} \leq 1$ respectively. The case of  $M(\z_p,1)$ was
settled in \cite{dydak-levin}.
\begin{theorem}{\rm ({\bf Dydak-Levin} \cite{dydak-levin})}
\label{moore-z-p}
A Moore space $M(\z_p, 1)$ is a classifying space
for finite dimensional compacta of $\dim_{\z_p} \leq 1$. 
\end{theorem}
The case of  $M(\z_{p^\infty}, 1)$ still remains open.
\begin{problem}
\label{moore-infinity}
Is a Moore space $M(\z_{p^\infty}, 1)$  a classifying space
for finite dimensional compacta of $\dim_{\z_{p^\infty}} \leq 1$?
\end{problem}

 Problem \ref{moore-infinity} leads  to:
\begin{definition}
A finite dimensional compactum $X$ is said to be extensionally regular
if for every prime $p$ we have
either $\dim_{\z_{p^\infty}} X \neq 1$ or  $\dim_{\z_{p^\infty}} X =1$
and $\edim X \leq M(\z_{p^\infty}, 1)$.
\end{definition} 

Note that every $M(G,1)$, $G\in \sigma$ is a $2$-dimensional
CW-complex and $M(G,n)=\Sigma^{n-1} M(G,1)$. Thus we obtain
\begin{corollary}
\label{corollary-moore}
For every $n$ and $G \in \sigma$ there  is an $(n-1)$-connected 
$(n+1)$-dimensional countable CW-complex  which is a classifying
space  for finite dimensional extensionall
regular compacta of $\dim_G \leq n$.

\end{corollary}

By a partial map of a space  $L$ to a  CW-complex $M$ we mean 
a map from a closed subset of $L$ to $M$. A collection ${\cal F}$ of partial maps
 from $L$ to $M$ is said to be representative if for every  closed 
subset $F'$ of $L$ and every map $f' : F' \lo M$, $M \in \cal M$
there is $f : F \lo M$  in $ {\cal F}$ such that $F' \subset F$ 
and
$f|F'$ is homotopic to $f'$. If $\cal M$ is a collection of CW-complexes
then 
a collection of partial maps of $L$ to the CW-complexes of $\cal M$
is said to be representative for $\cal M$ if it contains a representative collection
of partial maps to each $M \in \cal M$.
Note that 
 if $L$ is a compactum and $M$ is a countable CW-complex then there is
a countable representative collection ${\cal F}$ of partial maps  from $L$ to $M$
and a closed subset $A$ of $L$ is of $\edim \leq M$ if and only if
for every map $f : F \lo M$ in $ \cal F$, the map $f$ extends
over $A \cup F$.

In the proof of Theorem \ref{factorization} we will use the following
construction   from  \cite{l0,l1} for resolving partial maps.
     A  map between CW-complexes  is said to be  combinatorial if  
     the preimage of
   every subcomplex
   of the range is a subcomplex of the domain.
 Let $ L$ be  a simplicial complex and let  $ L^{[ m]}$        be
    the $ m$-skeleton of $ L$ (=the union of all simplexes of 
    $ L$ of $\dim \leq  m$).
By
 a resolution $EW( L, m)$   of $ L$   we mean a CW-complex $EW( L, m)$ 
 and
 a combinatorial map
 $\omega : EW( L, m) \lo   L$ such that $\omega$ is 1-to-1 over 
 $ L^{[ m]}$.
 Let $f : N \lo  M$  be a map of a subcomplex $N$ of $ L$ into
  a CW-complex $ M$.
A resolution $\omega : EW( L, m) \lo   L$
is said to resolve the map $f$   if
 the map  $f \circ\omega|_{\omega^{-1}(N)}$ extends
   to a map  $ f': EW( L, m) \lo  M$.  We will call $f'$  
   a  resolving map
   for $f$. The resolution is said to be
   suitable  for  a compactum $X$
if for   every simplex $\Delta$ of $ L$,
 $\edim X  \leq \omega^{-1}(\Delta)$.
     Note that if $\omega: EW( L, m) \lo  L$ is a resolution suitable
 for $X$   then  for every map  $\phi :  X \lo   L$  there is  a map   $\psi : X \lo EW( L, m)$
 such that  for every simplex $\Delta$ of $ L$ ,
  $(\omega \circ \psi)(\phi^{-1}(\Delta)) \subset \Delta$.
 We will call $\psi$ a combinatorial lifting of $\phi$.

Let $ L$   be a finite simplicial complex.
 Let       $f : N  \lo  M$ be a cellular  map  from a subcomplex $N$
of $ L$  to a CW-complex    $ M$ such that
$ L^{[ m]}\subset N$.
 Now   we will  construct  a
 resolution   $\omega: EW( L, m) \lo  L$  of
  $ L$  resolving    $f $ and 
   we will refer
 to this resolution as  the standard resolution for  $f$.
   We will associate     with the standard   resolution
 a  cellular resolving map $f' :  EW( L, m) \lo  M$ which will be
 called the standard resolving map.
  The standard  resolution is constructed  by induction
 on $l=\dim ( L\setminus N)$.

 For $ L=N$ set   $EW ( L, m)=   L$  and let    
 $\omega: EW( L, m) \lo  L$
 be the identity map   with    the standard resolving map $f'=f$.
 Let $l >  m$. Denote $ L'=N\cup  L^{[l-1]}$
 and   assume that     $ \omega': EW( L', m) \lo  L'$  is  the standard
 resolution     of       $ L' $
  for   $f$
  with the standard  resolving map  $f'  :  EW( L', m)  \lo  M$.
     The standard resolution      $\omega: EW( L, m) \lo  L$ 
     is constructed
 as follows.

 The  CW-complex $EW( L, m)$  is obtained  from  $EW( L', m)$  
 by attaching
   the mapping cylinder of  $f'|_{{\omega'}^{-1}({\partial \Delta})}$
 to      ${\omega'}^{-1}({\partial \Delta})$  for every  $l$-simplex $\Delta$ of $ L$
 which is not contained in $ L'$.
 Let $\omega : EW( L, m) \lo  L$ be the projection
which extends $\omega'$ by sending
  each mapping cylinder
 to  the corresponding $l$-simplex $\Delta$  such that
 the $ M$-part of the cylinder
 is sent to the barycenter  of     $\Delta$ and   each interval 
 connecting
 a point of ${\omega'}^{-1}({\partial \Delta})$
 with the corresponding point of the  $ M$-part of the cylinder is
 sent linearly to
 the interval connecting the corresponding point
 of ${\partial \Delta}$ with the barycenter of  $\Delta$.
 We can naturally  define  the extension of
 $f'|_{ {\omega'}^{-1}({\partial \Delta})}$ over its mapping cylinder
 by sending   each interval of the cylinder to the corresponding point 
 of $ M$.
 Thus        we define the  standard resolving map
 which extends $f'$
      over $EW( L, m)$.     The CW-structure of
   $EW( L, m)$  is induced  by the CW-structure of   $EW( L', m)$ and
           the natural CW-structures of
   the     mapping cylinders in  $EW( L, m)$.
  Then with respect to this CW-structure
    the  standard resolving map is cellular
   and $\omega$ is combinatorial.

It is easy to see  
from the construction of the standard resolution that
 for each simplex $\Delta$ of $ L$,
    $\omega^{-1}(\Delta)$ is either contractible or
  homotopy equivalent to $ M$ and
$\dim EW(L, m) \leq \dim L$ if $\dim M \leq  m+1$.
\begin{theorem}
\label{factorization-extensionally-regular}
 Any map from a finite dimensional extensionally regular
compactum $X$ 
to   a finite CW-complex $L$
 can be arbitrarily closely approximated by
a map  which factors
through a compactum $Z$ with $d_{Z} \leq d_X$ and $\dim Z \leq \dim L$.

\end{theorem}
{\bf Proof}.
Recall that a finite CW-complex is a compact ANR and hence
the identity map of a finite CW-complex can be arbitrarily closely
approximated by a map which factors through
a finite simplicial complex of the same dimension. Thus we may assume
that $L$ is a  finite simplicial complex. 

Let $g : X \lo L$ be a map. Set
$L_0 =L$ and $g_0 =g : X \lo L_0$. Fix $\epsilon>0$ and 
denote $l=\dim L$.
We will construct by induction
a sequence 
finite simplicial complexes $L_i$ with $\dim L_i \leq l$, bonding maps
$\omega^{i+1}_i : L_{i+1} \lo L_i$ and maps $g_i : X \lo L_i$
 such that $g_i$ and  $\omega^j_i \circ g_j $ are $\epsilon/2^j$-close
 for every $j>i$ 
where $\omega^j_i =L_j \lo L_i$ is the composition
 of the bonding  maps between $L_j, \dots,  L_i$  and
$\omega^i_i$ is the identity map of $L_i$. Denote 
 $Z=\invlim (L_i, \omega^{i+1}_i)$ and note that $\dim Z \leq l$.
 Also denote 
 $g'_i =\lim_{j \rightarrow \infty} \omega^j_i \circ g_j : X \lo L_i$ and
 note that $g'_i$ is a well-defined map and 
 $g'_{i+1} \circ g^{i+1}_i =g'_i$. Hence  the maps $g'_i$
 determine  the corresponding  map $g' : X \lo Z$ and for the projection
 $\omega_0 : Z \lo L_0=L$ we have that $g_0$ and $\omega_0 \circ g'$ are
 $\epsilon$-close. The construction will be carried out in such
 a way that 
 $d_Z \leq d_X$. Thus $Z$ and $\omega_0\circ g' $ 
 will  provide    the compactum and  the approximation 
 required in the theorem.
 
Assume that the construction is completed for $i$ and proceed to $i+1$ as follows.
Let   $m=\dim_G X$ for  a group $G \in \sigma$. Since the theorem is obvious for 
$\dim X =0$ we may assume that $m\geq 1$. By Corollary \ref{corollary-moore}
there is an $(m-1)$-connected $(m+1)$-dimensional  countable CW-complex $M$
classifying the finite dimensional
extensionally regular compacta of $\dim_G \leq m$.
Take a map $\alpha : F \lo M$ from a closed subset $F$ of $L_i$.
Replace the triangulation of $L_i$ by a sufficiently fine barycentric subdivision
so that 
$\alpha$ extends over a subcomplex $N$ of $L_i$ to
a map $f : N \lo M$  and 
for every simplex $\Delta$ of $L_i$ and every $j \geq i$
we have $\diam \omega^i_j (\Delta) \leq \epsilon/2^{i+1}$.
Since $M$ is $(m-1)$-connected we may assume that $N$ contains
the $m$-skeleton of $L_i$ and 
replacing $f$ by a cellular approximation  we also assume 
that
 $f : N \lo M$ is a cellular map. 
   Let $\omega : EW(L_i, m) \lo L_i$ be
the standard resolution resolving the map  $f$. Since $\edim X \leq M$
there is a combinatorial lifting  
$g_{i+1}  : X \lo EW(L_i,m)$ of $g_i$. Set $L_{i+1}$ to be
a finite subcomplex of $EW(L_i,m)$ containing $g_{i+1}(X)$ and
$\omega^{i+1}_i $ to be $\omega$  restricted
to $L_{i+1}$. Since the identity map
of $L_{i+1}$ can be arbitrarily closely approximated  
by a map which factors through a finite simplicial complex
of $\dim \leq \dim L_{i+1}$ we may assume that 
$L_{i+1}$ is a simplicial complex and the construction is completed.
Recall that  $\omega$ resolves the map $f$ and hence
 $\alpha \circ \omega^{i+1}_i|... : (\alpha \circ \omega^{i+1}_i)^{-1}(F) \lo M$
extends over $L_{i+1}$ as well.

Now we will show that the map $\alpha$ on the inductive
 step of the construction from $i$ to $i+1$ can
be chosen in a way that will  lead to $d_Z \leq d_X$.
Denote by ${\cal M}$ all the classifying  spaces mentioned
in Corollary \ref{corollary-moore} and  such that $\edim X \leq M$.
Note that $\cal M$ is a countable collection of countable CW-complexes.
Once $L_i$ is constructed take  
 a countable representative  collection ${\cal B}_i$ of
partial maps $\beta : B \lo M$ from closed subsets $B$ of $L_i$ to 
the CW-complexes $M$ of ${\cal M}$ and fix
a surjection $\tau_i : \N \lo {\cal B}_i $.
Take any bijection $\tau : \N \lo \N \times \N$
 such that for every $i$  and $\tau(i)=(j,k)$  we have $i\geq j$.
  Let $\tau(i)=(j,k)$ and $\beta=\tau_j^{-1}(k)\in {\cal B}_j$. Recall that $i \geq j$
and hence ${\cal B}_j$ is already constructed.
Thus
$\beta : B \lo M$  where $B$ is a closed subset of $L_j$
and $M \in \cal M$.
Denote 
$F=(\omega^{i}_j)^{-1}(B)$ and  set $\alpha=\beta \circ \omega_{i}^j : F \lo M$.
One can easily verify that choosing in such way the map $\alpha$ for constructing
$L_{i+1}$ leads to $\edim Z \leq M$ for every $M \in \cal M$ and hence
$d_Z \leq d_X$. The theorem is proved. $\black$
\\\\
In order to derive Theorem \ref{factorization} from
Theorem \ref{factorization-extensionally-regular} we need to by-pass
the difficulties imposed by Problem \ref{moore-infinity}.
We will need the extension versions of Theorems \ref{completion}
and \ref{decomposition}.
	\begin{theorem}
{\rm(\cite{Ol})}
\label{completion-extension} 
Let $K$ be a countable CW-complex and  
$X$ a space such that $\edim X \leq K$. Then
 there is a complete space $X'$ such that
$X \subset X'$ and $\edim X' \leq K$.
\end{theorem}

\begin{theorem}
{\rm (\cite{DrPacific})}
\label{decomposition-extension}
Let $K_1$ and $K_2$ be countable  CW-complexes, 
$K=K_1 * K_2$  the join of $K_1$ and $K_2$ and $X$ a compactum such
that $\edim X \leq K$. Then $X$ decomposes into $X=A \cup B$ such that
 $\edim A \leq K_1$ and $\edim B \leq K_2$.
\end{theorem}

In order to avoid a confusion with our previous use of the letter $L$
we will denote   the infinite dimensional lens space
model for  
 $K(\z_m,1)$ by ${\cal L}_m$ and, as ussual, 
 ${\cal L}^{[n]}_m$ stands for the
$n$-skeleton of ${\cal L}_m$ and $p$  for a prime number.
 By a projection from 
${\cal L}_{p^i}$ to ${\cal L}_{p^j}, j\geq i,$ we mean a cellular map
realizing 
the standard  monomorphism of  $\z_{p^i}$ into $\z_{p^j}$.
The restrictions of a projection to the skeletons of 
 ${\cal L}_{p^i}$ and ${\cal L}_{p^j}$ will be also called projections.
 Note that ${\cal L}^{[2]}_m =M(\z_m,1)$ and ${\cal L}^{[3]}_m$ is
 a $3$-dimensional lens space. 
 
 We assume that  $K(\z_{p^\infty},1)$ is represented
 as  the infinite telescope of a sequence of projections 
from ${\cal L}_{p^i}$ to ${\cal L}_{p^{i+1}}$
  and we consider  the lens spaces
	${\cal L}_{p^i}$ as subcomplexes of $K(\z_{p^\infty},1)$.
  
 It was shown in \cite{dydak-levin}
 that:
 \begin{theorem}{\rm (\cite{dydak-levin})}
 \label{dydak-levin-dim-3}
 Let $X$ be a finite dimensional 
 metric space  with $\dim_{\z_m} X \leq 2$ and
 $f : X \lo {\cal L}^{[n]}_m$  a map. Then there is a map 
 $f' : X \lo {\cal L}^{[3]}_m$ such that $f$ and $f'$ coincide on
 $f^{-1}({\cal L}^{[2]}_m)$.
 \end{theorem}

 \begin{proposition}
 \label{moore+}
Let $X$  be a  finite dimensional  compactum  with
 $\dim_{\z_{p^\infty}} X \leq 1$.  Then
 for every $i$ and every map $f: F \lo {\cal L}^{[2]}_{p^i}$
from a closed subset $F$ of $X$ there is $j \geq i$ such that
$f$ followed by a projection of ${\cal L}^{[2]}_{p^i}$ to ${\cal L}^{[2]}_{p^j}$
extends over $X$ as a map to ${\cal L}^{[3]}_{p^j}$.

 \end{proposition}
 {\bf Proof.} 
  Take a map $f : F \lo {\cal L}^{[2]}_{p^i}$ from
 a closed subset $F$ of $X$ and extend $f$ to 
 a map $h : X \lo K(\z_{p^\infty},1)$. Since $X$ is
compact, $h(X)$ is contained in a finite subtelescope of 
$K(\z_{p^\infty},1)$  and hence there is $j\geq i$ for which
$h$ can be homotoped to a map to $g : X \lo {\cal L}_{p^j}$
such that $g$ on $F$ coincides
with $f$ followed by a projection to ${\cal L}^{[2]}_{p^j}$.
Again by the  compactness of $X$
 there is $n$ such that $g(X) \subset {\cal L}^{[n]}_{p^j}$.
By Bockstein theorem and inequalities $\dim_{\z_{p^j}} X \leq 2$. 
Then,
 by Theorem \ref{dydak-levin-dim-3}, $g$ can be replaced by
a map to ${\cal L}^{[3]}_{p^j}$  which coinsides with $g$ on $F$
and the proposition follows.
$\black$
 \\
 \\
 Let $L$ be a simplicial complex.
 By the star  of a subset $A$  of 
 $L$ ,written $\st A$,  we mean the union of the simplexes of $L$
 which intersect $A$.  Let us say that for
 maps   $\phi : X \lo L$ and
 $\omega : Y \lo L$
  a map $\psi : X \lo Y$ is an almost combinatorial
 lifting of $ \phi $ to $Y$ if 
 $(\omega \circ \psi)(\phi^{-1}(\Delta)) \subset \st \Delta$.
 
 \begin{proposition}
 \label{main-moore-infinity}
 Let $X$ be a finite dimensional compactum with 
 $\dim_{\z_{p^\infty}} X \leq 1$, $L$ a finite simplicial complex with 
 $\dim L \leq 3$,
 $N$ a subcomplex of $L$, $\phi : X \lo L$  and
  $f : N \lo K(\z_{p^\infty},1)$  maps.
 Then there is a resolution $\omega : EW(L,1) \lo L$ such that
 $\dim EW(L,1) \leq 3$,
  $\omega$ resolves  the map $f$  and the map $\phi$ admits
  an almost combinatorial
 lifting  to $EW(L,1)$.
 \end{proposition}
 {\bf Proof.} Extending $f$ over the $1$-skeleton of $L$
 we  assume that $L^{[1]} \subset N$. 
 Since
 $f$ can be homotoped into ${\cal L}_{p^i} \subset K(\z_{p^\infty},1)$
  we may assume that $f(N) \subset {\cal L}_{p^i}$.
 By Bockstein Theorem and 
 Inequalities $\dim_{\z_{p^i}} X \leq 2$. Then,
 by Theorem \ref{dranishnikov-extension-theorem}, 
 $\edim X \leq  \Sigma M(\z_{p^i},1)=S^0 *  M(\z_{p^i},1)$ and
 hence, by Theorem \ref{decomposition-extension},
 $X$ decomposes into $X=A \cup B$ such that $\dim A \leq 0$
 and $\edim B \leq M(\z_{p^i},1)$, 
 and by Theorem \ref{completion-extension}
 we may assume that $B$ is $G_\delta$ and $A$ is $\sigma$-compact.

 Now replace the triangulation of $L$ by 
 its sufficiently fine subdivision.
Let 
$\cal R$ be the collection of the stars of the vertices of $L$
with respect to the barycentric subdivision $L_\beta$ of $L$.
Note that  $\cal R$ partions $L_\beta$ into
contractible subcomplexes with non-intersecting interiors.

By $\partial R, R \in \cal R,$ we denote the topological 
boundary of $R$ in $L$. 
 Clearly we may assume 
 that the triangulation of $L$ is so fine
that $f$ can be extended to  a  map  $f' : N' \lo {\cal L}_{p^i}$
over a subcomplex $N'$ of $L_\beta$ such that $N \subset N'$ and
 for every $R \in \cal R$
we have either   $R\subset N'$ or 
$(R \setminus \partial R) \cap N'=\emptyset$. Futhermore, we may assume that  $N'$
containes every  $1$-simplex  of $L_\beta$
 contained in $\partial R$ for every $R \in \cal R$.
Denote by $L'$ the subcomplex of $L_\beta$ which is the union of $N'$ with 
 $\partial R$ for all  $R \in \cal R$.
Let  $\cal R'$ be the collection of $R \in \cal R$
such that $R$ is contained in  $N'$, 
${\cal R}''={\cal R} \setminus \cal R'$,
$L''$ the subcomplex of $L'$ which is the union of
$\partial R$ for all $ R \in {\cal R}''$ 
and  $N''=N'\cap L''$.

Replace $f'$ by its cellular approximation  to ${\cal L}_{p^i}$. 
Then
$f'(N') \subset {\cal L}^{[3]}_{p^i}$ and 
  $f'(N'') \subset {\cal L}_{p^i}^{[2]}$. 
	Denote $f''=f'|_{N''} : N'' \lo {\cal L}_{p^i}^{[2]}$.
Let 
$\omega'': EW(L'',1) \lo L''$  be the standard 
resolution resolving the map $f''$  and let 
 $\rho'' : EW(L'',1) \lo {\cal L}^{[2]}_{p^i}$
be the standard  resolving  map for $f''$.
Extend $\omega''$ to the resolution
$\omega' : EW(L',1) \lo L'$ so that $\omega'$ is one-to-one over $N'$
and extend $\rho''$ to  the map 
$\rho' : EW(L',1) \lo {\cal L}^{[3]}_{p^i}$ defined by $f'\circ \omega'$ 
on $\omega'^{-1}(N')$. Then $\rho'$ resolves the map $f'$ over $L'$.

	Recall that the set $A$ is $\sigma$-compact and $0$-dimensional.
	Then one can  replace $\phi$ by 
 an  arbitrarily close approximation 
 and assume that for every $3$-simplex 
 $\Delta$ of $L$  and
 every $2$-simplex $\Delta'$ of the barycentric subdivision
 of $\Delta$ 
 we have that 
 $\phi(A) \cap (\Delta' \setminus \partial  \Delta)  =\emptyset$
 and hence 
 $\phi^{-1}(\Delta' \setminus \partial \Delta) \subset B$.
 Thus  for every
$2$-simplex 
$\Delta'$ of $L'$ we have that
$\phi^{-1}( \Delta'\setminus N') \subset B$  and
$\omega'^{-1}(\Delta')$ is either contractible or homotopy equivalent
to ${\cal L}^{[2]}_{p^i}$.
Then, since 
$\edim B \leq M(\z_{p^i},1)={\cal L}^{[2]}_{p^i}$, we conclude
that  $\phi$ restricted to
$X'=\phi^{-1}(L')$
admits a combinatorial
lifting $\psi' : X' \lo EW(L',1)$. 

Note that $(\rho'\circ \psi')(\phi^{-1}(\partial R)) \subset {\cal L}^{[2]}_{p^i}$
for $R \in {\cal R}''$.
 Then, by Proposition \ref{moore+}, for every $R \in {\cal R}''$ 
there is $j\geq i$ such that $\rho' \circ \psi'$ 
restricted to $\phi^{-1}(\partial R)$
and followed by a projection  $\tau : {\cal L}_{p^i}\lo {\cal L}_{p^j}$  
extends over $\phi^{-1}(R)$ as a map to 
${\cal L}^{[3]}_{p^j}$. Clearly $j$ can be replaced by any larger integer and
hence we can find  an integer $j$ that fits every $R \in  {\cal R}''$.
 Since $f'$ followed by $\tau$
 is homotopic to $f'$ as a map to $K(\z_{p^\infty}, 1)$
we can replace   $f'$ and $ \rho'$  by their compositions
with $\tau$ and assume that $f'$ and $\rho'$  are  maps to
 ${\cal L}^{[3]}_{p^j}$. 

Now define 
$EW(L,1)$ as the CW-complex obtained from $EW(L',1)$ by attaching
to $\omega'^{-1}(\partial R)$ 
the mapping cylinder of $\rho'$ restricted to $\omega'^{-1}(\partial R)$
 for every 
$R \in  {\cal R}''$.
Define $\omega : EW(L,1) \lo L$ as the map 
that sends  the ${\cal L}^{[3]}_{p^j}$-part of every mapping
cylinder to the vertex $v_R$  of $L$ contained in the corresponding 
set $R \in  {\cal R}''$ and linearly  extends $\omega'$ along
the intervals of the mapping cylinder and the corresponding intervals
connecting the points of $\partial R$ with $v_R$.
Let $\rho : EW(L,1) \lo {\cal L}^{[3]}_{p^j}$ be
a map naturally extending $\rho'$
over the mapping cylinders. Then $\dim E(L,1) \leq 3$,
$\omega$ is combinatorial,
$\rho$ resolves $f'$  and   hence $\rho$ resolves $f$ as well.

Finally note that  for every $R \in {\cal R}''$ the map $\psi'$
restricted to  $\phi^{-1}(\partial R)$  and considered as a map
to  the mapping cylinder attached to $\omega'^{-1}(\partial R)$
can be homotoped to 
the ${\cal L}^{[3]}_{p^j}$-part of  the mapping cylinder and then extended over 
$\phi^{-1}(R)$ as a map to ${\cal L}^{[3]}_{p^j}$. This way we can extend
$\psi'$  to  a map  $\psi : X \lo EW(L,1)$ so that  
$(\omega \circ\psi) (\phi^{-1}(R) ) \subset R$ for every
$R \in \cal R$.

Thus with respect to the original triangulation of $L$, the map
$\omega: EW(L,1) \lo L$ is a resolution resolving the map $f$
and admitting 
  an almost
combinatorial lifting $\psi : X \lo EW(L,1)$ of $\phi$,  
 and the  proposition follows.
$\black$
\\
\\
{\bf Proof of Theorem \ref{factorization}}. The proof of 
Theorem \ref{factorization-extensionally-regular} applies to
prove Theorem \ref{factorization} with the following  adjustment: 
 for every prime $p$  
with $\dim_{\z_{p^\infty}} X =1$ replace  the Moore space $M(\z_{p^\infty},1)$
by the Eilenberg-Mac Lane complex  $K(\z_{p^\infty},1)$ and
the standard resolution of partial 
maps to  $K(\z_{p^\infty},1)$ by the the resolution from
Proposition \ref{main-moore-infinity}. $\black$
 
\end{section}
\begin{section}{Main result}
\label{main-result}
The goal of this section is to prove Theorem \ref{main}. The following results
will be used in the proof.

Using Alexander duality and Kunneth formula in the Leray form one can show:
\begin{theorem}
{\rm (\cite{dranishnikov-sbornik})}
\label{duality} 
Let $X$ and $Y$ be compacta such that  $Y\subset \re^n$  and 
$\dim X \times Y \leq n-1$. Then for every open ball $U$ in $\re^n$
and $i \geq 0$ 
we have $\dim_{H_i(U\setminus Y) } X \leq i$ where $H_*(U\setminus Y)$ is
the reduced integral homology.
\end{theorem}
Let us remind that 
a tame compactum $X \subset \re^n$ with $\dim X \leq n-3$ is characterized by 
the following property: 
for every open ball $U$ in $\re^n$ the complement  $U \setminus X$ is simply connected.
Using Theorems \ref{duality} and \ref{dranishnikov-extension-theorem} one can show:
\begin{theorem}
\label{dranish-tame}
{\rm (\cite{dranishnikov-sbornik})}
Let   $X$ be a finite dimensional compactum  and 
 $ Y \subset \re^n$      a tame compactum of  $\dim Y\leq n-3$
such  that $\dim X \times Y \leq n-1$.
Then every map $g : X\lo \re^n$ can be arbitrarily closely approximated
by a map $g' : X \lo \re^n$ so that $g'(X) \cap Y =\emptyset$.
\end{theorem}

\begin{proposition}
\label{tame-dim-3}
Let $X$ be a compactum with $\dim X \leq 3$ and $n\geq 6$. Then 
any map from $ X$ to $\re^n$ can be arbitrarily closely approximated
by 
 a map $f : X \lo \re^n$ so that $d_{f(X)} \leq d_X$ and
$f(X)$ is tame in $\re^n$.
\end{proposition}
{\bf Proof.} By Stanko's re-embedding theorem \cite{stanko} it  suffices
to construct $f$ with $d_{f(X)} \leq d_X$.
 The cases 
 $n>6$  or $\dim X \leq 2$ are
 trivial since any map from $X$ to $\re^n$ can be approximated
by an embedding. So the only case  we need to consider  is $\dim X =3$ and
$n=6$.

Let $g : X \lo \re^n$ be a map and let $ \alpha_F : F \lo K(G,m), m\geq 1,$ be a map
from a compact subset $F$ of $\re^n$ to  $K(G,m)$ such that
$\dim_G X \leq m$ and $G$ is a group in the Bockstein basis $\sigma$.
Extend $ \alpha_F$ over a closed neighborhood $F^+$ of $F$
to a map $ \alpha_F^+ : F^+ \lo  K(G,m)$ and
extend $g$ restricted to $g^{-1}(F^+ )$ and followed by $ \alpha_F^+$
to  a map $\beta : X \lo K(G,m)$. 
Note that, since $X$ is compact, the map $\beta$ can be considered as a map to
a finite subcomplex $K$ of $K(G,m)$. 
Let $\epsilon>0$ be such that any two $2\epsilon$-close maps 
to $K$ are homotopic.
Approximate $g$  through
 a $3$-dimensional finite simplicial complex $L$ 
and  maps $\gamma : X \lo L$ and $g_L : L \lo \re^n$ such that
$\gamma$ is surjective and:

(i) $g$ and $g_L \circ \gamma$ 
are so close that
 $g(\gamma^{-1}(F_L)) \subset F^+$ for $F_L =g_L^{-1}(F)$ and
the maps $ \alpha_F^+ \circ g$ and
$ \alpha_F^+ \circ g_L \circ \gamma$ restricted to $\gamma^{-1}(F_L)$
are $  \epsilon$-close;

(ii) the fibers of $\gamma$ are so small that
that there is a map $\beta_L : L \lo K$ so that
$\beta$ and $\beta_L \circ \gamma$ are $  \epsilon$-close.

Take $y\in F_L$ and let $x\in X$ be such that $\gamma(x)=y$.
By (ii), $\beta(x)$ and $(\beta_L\circ \gamma)(x)=\beta_L(y)$ 
are $  \epsilon$-close.
By (i), $\beta(x)=(\alpha^+_{F} \circ g)(x)$ and
$(\alpha^+_F \circ g_L \circ \gamma) (x)=
(\alpha^+_F \circ g_L)(y)=(\alpha_F \circ g_L)(y) $  
are $  \epsilon$-close.
Thus we get that $\beta_L$ and $\alpha_{F} \circ g_L$
are $2\epsilon$-close on $F_L$ and 
 hence  replacing $\beta_L$ by a homotopic map
we may assume that $\beta_L$ coincides  on $F_L$ with $g_L$ followed by $\alpha_F$.

Since $\dim L \leq 3$  and $n\geq 6$ we can in addition assume 
that $g_L$ is   finite-to-one and  $g_L$ is
not one-to-one over  only finitely many points of $g_L(L) $.
 Then, since $K(G,m)$ is a connected CW-complex,
we can change (up to homotopy)  $\beta_L$  outside the set  $F_L$
so that $\beta_L$ will be constant on each fiber of $g_L$ and  hence we
can assume 
that $\beta_L$ factors through a map $\alpha : g_L(L) \lo K(G,m)$
so that $\alpha$ coincides with $\alpha_F$ on $g_L(L) \cap F$.
Thus  the map $g$ can be arbitrarily closely approximated 
 by a map $f =g_L \circ \gamma : X \lo \re^n$
so that   $\alpha_F$ extends  over $F \cup f(X)=F \cup g_L(L)$.  
 Note that there is  
a countable representative collection of partial maps $\alpha_F$ 
from compact subsets of $\re^m$ to Eilenberg-MacLane complexes $K(G,m)$
with $\dim_G X \leq m$, $G \in \sigma$.  Then the above property implies
that $g$ can be arbitrarily closely approximated by a map $f: X \lo \re^n$
with $d_{f(X)} \leq d_X$  and the proposition follows.
$\black$
\\\\
 Note that the use of Stanko's re-embedding theorem in Theorem \ref{tame-dim-3}
can be easily avoided 
 by constructing $f(X)$ to be the intersection
of a decreasing sequence  of sufficiently close PL-regular neighborhoods  of 
$3$-dimensional  finite simplicial complexes in $\re^n$.

\begin{proposition}
\label{plus-construction}
Let $K$ be a CW-complex. Then one can attach to $K$  cells of 
dimensions $\leq 3$
to obtain a simply connected CW-complex $K^+$  so that the inclusion
of $K$ into $K^+$ induces an isomorphism of the integral homology in
dimensions $ >1$.
\end{proposition}
{\bf Proof.} 
Clearly by attaching intervals to $K$ we can turn $K$ into
a connected CW-complex preserving the integral homology of 
$K$ in dimensions $>0$. Thus we may assume that $K$ is connected.

Let a simply connected
CW-complex $K'$ be obtained from $K$ by attaching $2$-cells
to kill the fundamental group of $K$. Consider 
the inclusion $i: K \lo K'$ and 
the quotient map
$p : K' \lo K'/K$.
Note that
$i_*(H_2(K)) =\ker p_*$ 
and  $H_2(K'/K)$ is a free group
since  $K'/K$ is   a bouquet of $2$-spheres. Then $p_*(H_2(K'))$ is a free group as well.
Take a collection $\omega_j, j \in J$ of $2$-cycles of $K'$ such that
$p_*[\omega_j], j \in J$ are free generators of $p_*(H_2(K'))$.
Since $K'$ is simply connected one can enlarge $K'$ to a CW-complex $K''$
by attaching for every $j$ a $3$-cell  $\Omega_j$  such that $\partial \Omega_j$
is homologous in $K'$ to $\omega_j$
(we consider the cellular homology). Replacing $\omega_j$ by $\partial \Omega_j$
 assume that $\omega_j=\partial \Omega_j$. Let us show
that the inclusion of $K$ into $K''$ induces an isomorphism of the homology groups
in dimensions$>1$.

Take a $2$-cycle $\alpha$ of $K''$. Then $\alpha$  lies in 
$K'$ and   $p_*[\alpha]=\sum n_j p_*[\omega_j], n_j \in\z$. Thus $\alpha -\sum n_j \omega_j$
is homologous to a cycle in $K$ and since $\omega_j=\partial \Omega_j$
we get that $\alpha$ is homologous in $K''$ to a cycle in $K$. Hence
the inclusion of $K$ into $K''$ induces an epimorphism of $2$-homology.

Take a $2$-cycle $\alpha$ in $K$  homologous to $0$ in $K''$. Then
$\alpha=\partial (\beta +\sum n_j \Omega_j)$  with $\beta$ being a $3$-chain in $K$.
Thus $\alpha= \partial \beta  +\sum n_j \omega_j $ and
$p_*[\alpha]=\sum n_j p_*[\omega_j]=0$ and hence $n_j=0$ for every $j$.
Thus $\alpha =\partial \beta$ is homologous to $0$ in $K$ and hence
the inclusion of $K$ into $K''$ induces a monomorphism of $2$-homology.

Since we attached to $K$ only cells of $\dim \leq 3$ the inclusion
of $K$ into $K''$ induces a monomorphism of $3$-homology. Take
a $3$-cycle $\alpha$ in $K''$. Then $\alpha =\beta +\sum n_j \Omega_j$
where $\beta$ is a $3$-chain in $K$. Hence
$p_*[\partial \alpha ]=p_*[\partial \beta +\sum n_j \omega_j]=
\sum n_j p_*[\omega_j]=0$ and we get that $n_j=0$ for every $j$. 
Thus $\alpha=\beta$ and hence the inclusion of $K$ into $K''$ 
induces an epimorphism of $3$-homology.

Clearly the homology groups of $K$ and $K''$ coincide in $\dim >3$.
Thus  the proposition holds with $K^+=K''$.
$\black$

\begin{proposition}
\label{type-approximation}
Let $X$ be a  compactum with $\dim X \leq n-2$ and 
$n \geq 6$. Then any map $f: X \lo \re^n$ can be arbitrarily closely
approximated by a map $f' : X \lo \re^n$ so that $d_{f'(X)} \leq d_X$.
\end{proposition}
{\bf Proof.} We prove the proposition by induction on $\dim X$.
 Clearly the proposition holds $\dim X=0$. Assume that
the proposition is proved for compacta of $\dim \leq \dim X -1$.
Fix $\epsilon >0$ and
partition $X$ into finitely many closed subsets $X=\cup X_j$ with disjoint
interiors such that $\dim \partial X_j \leq \dim X -1$ and
$f(X_j)$ is contained in an open $\epsilon$-ball in $\re^n$.
By the induction hypothesis  we can assume that $d_{f(\partial X_j)} \leq
d_{X_j} \leq d_X$. 

Fix $j$ and  an open $\epsilon$-ball $U$  containing  $X_j$.
Take a map $\alpha_F : F \lo K(G,m)$  from a compact subset $F$ of $U$
to an Eilenberg-MacLane
complex  $K(G,m)$ such that $\dim_G X\leq m$, $G \in \sigma$.
 By Proposition \ref{type-minus} and Theorem \ref{decomposition}
decompose $\re^n$ into $\re^n=A \cup B$
with $d_A \leq d_X $ and $d_B \leq n-1 \ominus d_X$. Clearly we may assume that
$f(\partial X_j) \subset A$. Extend $\alpha_F$ to a map $\alpha_W : W \lo K(G,m)$
over an open subset $W$ of $\re^n$ such that $A \subset W$ and
  let $Y$ be the closure of $U \setminus W$ in $\re^n$.  
	Then $Y $ is a compact  subset  of $B$
and, 
 by Theorem \ref{duality} and Proposition \ref{type-minus}, for the set
$V=U\setminus Y=U\cap W$ we have 
$\dim_{H_i(V)} X \leq i$ for every $i$.  Take a triangulation 
of $U$ and, Proposition \ref{plus-construction}, attach to $V$  
cells of dimensions
 $\leq 3$ to obtain a simply connected CW-complex $V^+$ 
preserving the homology of $V$ in dimensions$>1$.  Then, 
by Theorem \ref{dranishnikov-extension-theorem},
$\edim X \leq V^+$ and hence here is a map $\phi: X_j \lo V^+$
extending $f$ restricted to $\partial X_j$. 

Let $L$ be a finite subcomplex 
of $V^+$ containing $\phi(X_j)$. Denote $L^-=V \cap L$,
$L^+$=the $3$-skeleton of $K$, $X^-_j=(\phi)^{-1}(L^-)$
and $X^+_j =(\phi)^{-1}(L^+)$.
 Clearly $L^-$ and  $L^+$ are subcomplexes of $L$,
$\dim L^+\leq 3$ and $L=L^-\cup L^+$. 
Take any map $\psi : L \lo U$ such that $\psi$ does not move the points
of $L^-$and take an open neighborhood $V^-$ of $L^-$ in $V$ so that 
the closure of $V^-$ is contained in $V$.
By Theorem \ref{factorization}
approximate $\phi$ restricted to $X_j^+$ by the composition
of  maps $\phi^+_Z :  X^+_j \lo Z $ and $\phi^+_L : Z \lo L^+$ such that
 $\edim Z \leq d_X$ and $\dim Z \leq 3$.
By Proposition \ref{tame-dim-3},  the map $\phi^+_L$ followed by $\psi$
can be approximated by a map to $U$ with the image of 
dimension type $\leq d_Z\leq d_X$. This trivially implies that
$f$ restricted to $X_j$ can be
arbitrarily closely approximated  by a map $f'_j : X_j \lo U$ such that
$f'_j$ coincides on $\partial X_j$ with $f$ and $f'_j(X_j) \setminus V^-$ is
of dimension type $\leq d_X$. Then $\alpha_W$ restricted to the closure of
$V^-$ extends over $ f'_j(X_j) \setminus V^-$. Thus we can extend
$\alpha_F$ to a map $\alpha : F \cup f'_j(X_j) \lo K(G,m)$. Since
there is a countable representative collection of  maps
from compact subsets of $\re^n$ to the Eilenberg-MacLane complexes
for the groups in $\sigma$
 the last property
allows 
 one to  achieve that $f'_j(X_j)$ is of dimension type $\leq d_X$.
Assemble the maps $f'_j$ into a map $f ' : X \lo \re^n$ which 
provides the required $\epsilon$-approximation of $f$ and  the proposition follows.
$\black$

\begin{theorem}
\label{main-theorem}
Let $X$ and $Y$ be  
 compacta such that $\dim X \leq n-2$, $\dim Y \leq n-2$,
$\dim X \times Y \leq n-1$ and $n\geq 6$. Then $X$ and $Y$ do not
admit a stable intersection in $\re^n$.
\end{theorem}
{\bf Proof.}
We prove the theorem by induction on $\dim X$. The case $\dim X=0$ is trivial.
Assume that $\dim X >0$. Take maps $f : X \lo \re^n$ and $g : Y \lo \re^n$.
Fix $\epsilon >0$ and
partition $X$ into finitely many closed subsets $X=\cup X_j$ with disjoint
interiors such that $\dim \partial X_j \leq \dim X -1$ and
$f(X_j)$ is contained in an open $\epsilon$-ball in $\re^n$.
By the induction hypothesis  and Proposition \ref{type-approximation}
 we can  replace $f$ and $g$ by 
arbitrarily close approximations and 
assume that $d_{g(Y)} \leq d_Y$ and 
 $f(\partial X_j) \cap g(Y)= \emptyset$  for every $j$.
Fix $j$ and  an open $\epsilon$-ball $U$  containing  $X_j$.
By Theorem \ref{duality} and Proposition \ref{type-minus}, for the set
$V=U\setminus g(Y)$ we have 
$\dim_{H_i(V)} X \leq i$ for every $i$.  Take a triangulation 
of $U$ and, Proposition \ref{plus-construction}, attach to $V$  
cells of dimensions
$\leq 3$ to obtain a simply connected CW-complex $V^+$ 
preserving the homology of $V$ in dimensions$>1$.  Then, 
by Theorem \ref{dranishnikov-extension-theorem},
$\edim X \leq V^+$ and hence here is a map $\phi: X_j \lo V^+$
extending $f$ restricted to $\partial X_j$. 

Let $L$ be a finite subcomplex 
of $V^+$ containing $\phi(X_j)$. Denote $L^-=V \cap L$,
$L^+$=the $3$-skeleton of $K$, $X^-_j=(\phi)^{-1}(L^-)$
and $X^+_j =(\phi)^{-1}(L^+)$.
 Clearly $L^-$ and  $L^+$ are subcomplexes of $L$,
$\dim L^+\leq 3$ and $L=L^-\cup L^+$. 
Take any map $\psi : L \lo U$ such that $\psi$ does not move the points
of $L^-$and take an open neighborhood $V^-$ of $L^-$ in $V$ so that 
the closure of $V^-$ is contained in $V$.
By Theorem \ref{factorization}
approximate $\phi$ restricted to $X_j^+$ by the composition
of  maps $\phi^+_Z :  X^+_j \lo Z $ and $\phi^+_L : Z \lo L^+$ such that
 $\edim Z \leq d_X$ and $\dim Z \leq 3$.
By Proposition \ref{tame-dim-3},  the map $\phi^+_L$ followed by $\psi$
can be approximated by a map to $U$ with a $3$-dimensional tame image of 
dimension type $\leq d_Z\leq d_X$. This trivially implies that
$f$ restricted to $X_j$ can be
arbitrarily closely approximated  by a map $f'_j : X_j \lo U$ such that
$f'_j$ coincides on $\partial X_j$ with $f$ and $f'_j(X_j) \setminus V^-$ is
a tame $3$-dimensional compactum 
of dimension type $\leq d_X$. Then, by Proposition  \ref{tame-dim-3},
$g$ can be arbitrarily closely approximated by a map $g' : Y \lo \re^n$
such that $(f'(X_j)\setminus V^- )\cap g'(Y)  =\emptyset$.  Recall that
$g(Y)$ does not meet the closure of $V^-$  and hence we can assume 
that $f'(X_j) \cap g'(Y) =\emptyset$.

The above procedure  can be carried out consecutively for all $X_j$ to obtain
approximations of $f$ and $g$ with disjoint images
and
the theorem follows.
$\black$
\\\\
{\bf Proof of Theorem \ref{main}.}
 Theorem \ref{main} follows from Theorems \ref{main-theorem},
\ref{stable} and \ref{complementary}. $\black$

\end{section}
\begin{section}{ Remarks}
Problem \ref{moore-infinity} can be considered
in a more general context.

\begin{problem}
\label{general-extension-problem}
Let $K$ be a connected CW-complex whose fundamental group is
abelin and $X$ a finite dimensional compactum 
such that $\dim_{ H_n(K)} X \leq n$ for every $n>0$.
Does it imply that $\edim X \leq K$?
\end{problem} 
It turns out that this problem reduces to $3$-dimensional 
compacta.

\begin{proposition}
\label{geeneral-extension-to-3}
Problem \ref{general-extension-problem} is equivalent to the same
problem with the additional assumtion that $\dim X \leq 3$.
\end{proposition}
{\bf Proof.} 
Let $f : F \lo K$ be a map form a closed 
subset $F$ of $X$. By Proposition \ref{plus-construction}, 
$K$ can be enlarged to a simply connected CW-complex $K^+$
preserving the integral homology of $K$ in dimensions $>1$ 
with $\dim K^+ \setminus K  \leq 3$. Then, 
by Theorem \ref{dranishnikov-extension-theorem},
$f$ extends to a map $g : X \lo K^+$. Let $L$ be the $3$-skeleton
of $K^+$ and $Y=g^{-1}(L)$.  By, Theorem \ref{factorization},
$g$ restricted to $Y$ factors up to homotopy through 
a compactum $Z$  with $\dim Z \leq 3$ and $d_Z \leq d_X$
and
a map $\psi : Z \lo L$. Now assuming  
that $\psi$ restricted to $\psi^{-1}(K \cap L)$ extends over $Z$
as a map to $K$ we get that $f$ extends over $X$ as a map to $K$
 as well and the proposition follows. $\black$
 \\
 \\
 Let us state without a proof  a few more results
 related to the techniques presented in this paper.
 \begin{itemize}
 \item Proposition \ref{factorization} can be extended to the following
 result. Let $X$ be a finite dimensional dimensional compactum,
  $L$  a finite CW-complex and $f : X \lo L$
 a map. Then $f$ can be arbirarily closely approximated by 
 a map that factors through a compactum $Z$ with
 $\dim Z \leq \max \{ \dim L ,3 \}$ and $d_Z \leq d_X$.
 
 \item
 It turns out that Theorem \ref{factorization-extensionally-regular}
 and Problem \ref{moore-infinity} are closely related. Namely,
 a Moore space $M(\z_{p^\infty},1)$ is a classifying space
 for finite dimensional compacta  with
 $\dim_{\z_{p^\infty}} = 1$ if and only
 if for every $3$-dimensional compactum $X$
 with $\dim_{\z_{p^\infty}} X=1$ and every map
 $f : X \lo \re^2$ we have that $f$ can be arbitrarily closely
 approximated by a map that factors through 
 a compactum $Z$ with $\dim Z \leq 2$ and $\dim_{\z_{p^\infty}} Z\leq 1$.
 
 \item
 Theorem \ref{factorization-extensionally-regular} admits the following
 generalization.
 Let $f: X \lo Y$ be a map of finite dimensional compacta
 such that $X$ is extensionally regular.
 Then $f$ factors through a compactum $Z$ with
 $\dim Z \leq \dim Y$ and $d_Z \leq d_X$.
 
 \end{itemize}

\end{section}

Michael Levin\\
Department of Mathematics\\
Ben Gurion University of the Negev\\
P.O.B. 653\\
Be'er Sheva 84105, ISRAEL  \\
 mlevine@math.bgu.ac.il\\\\

\begin{thebibliography}{99}


\bibitem{dranish-west}
A. Dranishnikov and J. West, 
{\em On compacta that intersect unstably in Euclidean space}, Topol.
Appl. 43 (1992), 181-187.


\bibitem{DrPacific} Dranishnikov A. N. {\em On the mapping intersection problem.} Pacific J. Math. 173 (1996), 403-412.
\bibitem{dranishnikov-intersection}
   Dranishnikov, A. N. {\em On the dimension of the product of two compacta
   and the dimension of their intersection in general position in Euclidean space.}
   Trans. Amer. Math. Soc.  352  (2000),  no. 12, 5599--5618.

\bibitem{d-r-s}
A. Dranishnikov, D. Repovs, E. Scepin, {\em  On intersection of compacta 
in Euclidean space: 
the metastable case}, Tsukuba J. Math. 17 (1993), 549-564.


\bibitem{DrArxiv} Dranishnikov, A. N. {\em Cohomological dimension theory of compact 
metric spaces,} Topology
Atlas invited contribution, http://at.yorku.ca/topology.taic.html 
(see also  arXiv:math/0501523).
 \bibitem{first-exotic} Dranishnikov A., Levin M. 
{\em Dimension of the product and classical formulae of dimension theory,}
 Trans. Amer. Math. Soc., to appear.


\bibitem{DrUsp} Dranishnikov, A. N. {\em Homological dimension theory.} Russian Math. 
Surveys 43 (4) (1988), 11-63.
\bibitem{dranishnikov-sbornik}
A.N. Dranishnikov {\em An extension of mappings into CW-complexes}, 
Mat. Sb. 182 (1991), 1300-1310; English transl., Math. USSR Sb. 74 (1993), 47-56.

\bibitem{Dy}
 Dydak, Jerzy {\em Cohomological dimension and metrizable spaces. II.} 
Trans. Amer. Math. Soc. 348 (1996), no. 4, 1647--1661.

\bibitem{DyA}
Dydak, Jerzy {\em Algebra of dimension theory.} Trans. Amer. Math. Soc. 358 (2006), 
no. 4, 1537--1561.



\bibitem{dydak-levin}
J. Dydak and M. Levin, {\em Extensions of maps to $M(\z_m, 1)$}, preprint (2013),
arXiv:1310.1711

\bibitem{Ku}
Kuzminov, V. I. {\em Homological dimension theory.} Russian Math Surveys 23 (5)  (1968), 
1-45.

\bibitem{m-rubin}
D. McCullough and L. Rubin, {\rm Some m-dimensional compacta admitting a dense set of 
imbedding
into $R^{2m}$}, Fund. Math. vol 133 (1989), 237-245.

\bibitem{l0}Michael Levin, {\em Constructing compacta of different extensional 
dimensions}.
Canad. Math. Bull. 44(2001), no. 1, 80–86.
\bibitem{l1}
Michael Levin, {\em Acyclic resolutions for arbitrary groups}. Israel J. Math. 135
(2003), 193–203.

\bibitem{Ol} Olszewski W. {\em Completion theorem for cohomological dimensions.} 
Proc. Amer. Math. Soc. 123 (1995) 2261--2264.
\bibitem{Sch}
Shchepin, E. V. {\em Arithmetic of dimension theory.} Russian Math. Surveys 53 (1998), 
no. 5, 975--1069.

\bibitem{stanko} M.A. Stanko,
 {\em Embeddings of compacta in Euclidean space}, Math. USSR-Sbornik 10 (1970),
234-254.
\bibitem{t-s}
S. Spiez and H. Torunczyk,
 {\em Moving compacta in $R^n$ apart}, Topol. Appl. 41 (1991), 193-204.










\end{thebibliography}
\end{document}